\documentstyle[12pt,amsfonts]{article}
\title{Topological rigidity of automorphism actions on nilmanifolds}
\author{Siddhartha Bhattacharya}
\begin{document}
\date{}
\maketitle
\baselineskip=20pt
\section{Introduction}

Let $G$ be connected simply connected nilpotent Lie group and $D$
be discrete uniform subgroup of $G$. Then 
$X = G/D $ is called a nilmanifold. 
If $X_{1} = G_{1}/D_{1}$ and
$X_{2} = G_{2}/D_{2}$ are nilmanifolds then a map 
$f:X_{1}\rightarrow X_{2}$ is said to be a ${\it homomorphism \/}$ if it
is induced by a continuous
homomorphism from $G_{1}$ to $G_{2}$ which maps $D_{1}$ into $D_{2}$.
Isomorphisms and automorphisms are defined similarly.
A map $f:X_{1}\rightarrow X_{2}$ is said to be ${\it affine \/}$ if there
exists an element $a$ in $G_{2}$ and a continuous homomorphism
$A :G_{1}\rightarrow G_{2}$ such that 
$f(gD_{1}) = aA(g)D_{2}$ for all $g$ in $G_{1}$. If $\Gamma$ is a discrete
group and $X$ is a nilmanifold then a $\Gamma$-action $\rho $ on $X$ is
said
to be an ${\it automorphism\ action \/}$ if each $\rho(\gamma )$ is an
automorphism of $X$.

It is known that any nilmanifold $X = G/D $ carries a unique
$G$-invariant probability measure $\lambda_{X}$.
It is easy to see that for any nilmanifold $X$, the measure $\lambda_{X}$ is
invariant under any affine
action on $X$. An automorphism action $\rho$ of a discrete group
$\Gamma$ on a nilmanifold $X$ is said to be ${\it ergodic \/}$ if for
every $\Gamma$-invariant function $f$ in $L^{2}(X,\lambda_{X})$ is a
constant almost everywhere.

Suppose  $X_{1}, X_{2}$ are nilmanifolds
and $\rho ,\sigma$ are continuous actions of a discrete group $\Gamma$ on
$X_{1}$ and $X_{2}$ respectively. 
If  $f:X_{1}\rightarrow X_{2}$ is a continuous map and $\Gamma_{0}$ is 
a subgroup of $\Gamma$ then $f$ is said to be
{\it $\Gamma_{0}$-equivariant \/} if  
$ f\circ \rho (\gamma ) = \sigma (\gamma )\circ f, \ \ \forall
\gamma
\in{\Gamma_{0} }$.
A continuous map $f:X_{1}\rightarrow X_{2}$ is said to be
{\it almost equivariant \/} if
 there exists a finite-index subgroup $\Gamma_{0}\subset\Gamma $ such
that $ f$ is $\Gamma_{0}$-equivariant.

In this paper we show that if $\rho ,\sigma $ are automorphism actions
 on nilmanifolds satisfying certain conditions
then every $\Gamma$-equivariant continuous map from $(X_{1},\rho )$ to
$(X_{2},\sigma )$ is an affine map.
 When $\Gamma = {\mathbb{Z}}$,
$(X_{2},\sigma)$ is a factor of $(X_{1},\rho)$ and $\rho, \sigma $ 
are generated by affine transformations, this
phenomenon has
been studied in
[AP],[Wa1] and
[Wa2]. In this case  a necessary and sufficient condition for 
existence of a non-affine  $\Gamma$-equivariant map is given in 
[Wa2]. Our methods are however different and applicable in more
general situations.

This paper is organized as follows.
In section 2 we study structure of continuous equivariant maps from 
 $(X_{1},\rho )$ to $(X_{2},\sigma )$. In Theorem 1 we give a necessary
and
sufficient condition for
 existence of a non-affine almost equivariant map.
 For any Lie group $G$, by $L(G)$ we denote the Lie algebra of $G$.
If $\rho$ is an automorphism action of a discrete group $\Gamma$ on a
 nilmanifold $X = G/D$, then by $\rho_{e}$ we denote the 
$\Gamma$-action on $L(G)$ induced by $\rho$. We prove the following.
\newline
\newline
{\bf Theorem 1 : }
 {\it Let  $X_{1} = G_{1}/D_{1}, X_{2} = G_{2}/D_{2}$ be nilmanifolds and
$\rho ,\sigma $
be automorphism actions of a discrete group $\Gamma $  on
 $X_{1}$ and $X_{2}$ respectively. Then
there exists a non-affine almost equivariant continuous  map
from
$(X_{1},\rho )$ to $(X_{2},\sigma )$ if and only if
the following two conditions are satisfied.
\newline 

a) There exists a non-constant $ \Gamma$-invariant continuous
function from 

$\ \ (X_{1}, \rho)$ to $\mathbb{R}$.

b) There exists a nonzero vector $v$ in $L(G_{2})$ with
 finite $\sigma_{e}$-orbit.
\/}
\newline
\newline 
As a consequence we obtain that if either $(X_{1},\rho )$ is ergodic
or $(X_{2},\sigma )$ is expansive then every
continuous $\Gamma $-equivariant map
from
$(X_{1},\rho )$ to $(X_{2},\sigma )$ is an affine map (see corollary
2.1). 

In section 3  we  consider the case when $X_{1}$ is a torus.
If 
$\rho $ is an automorphism action of a discrete group $\Gamma $ on
a torus $T^{m}$ then we denote
 the
 induced automorphism action  of $\Gamma $ on the dual group 
 ${\widehat{T}}^{m}$  by 
$\widehat{\rho}$ . By $F_{\rho}$ we denote the
subgroup of
${\widehat{T}}^{m}$
which consists of all elements with finite $\widehat{\rho}$-orbit
 and by $\Gamma_{\rho}$ we denote the subgroup of $\Gamma$ 
consisting of all elements which acts trivially on $F_{\rho}$ under
the action $\widehat{\rho}$. 
We prove the following.
\newline
\newline
{\bf Theorem 2  : } {\it Let $\Gamma $ be a discrete group, 
$T^{m}$ be the $m$-torus and $X = G/D$ be a nilmanifold. Let 
$\rho ,\sigma $
be automorphism actions of $\Gamma $  on
$T^{m}$ and $X$ respectively.  Then
there exists a non-affine continuous $\Gamma $-equivariant map
from
$(T^{m},\rho )$ to $(X,\sigma )$ if and only if
the following two conditions are satisfied.
\newline

a)  $\ (T^{m},\rho )$ is not ergodic.

b) There exists a nonzero vector $v$ in $L(G)$ which is fixed
  by $\Gamma_{\rho}$ under 

$\ \ \ $ the action $\sigma_{e}$.
 \/}
\newline
\newline
In section 4 we consider the case when $\Gamma$ is abelian and 
$(X_{2},\sigma)$ is a topological factor of $(X_{1},\rho)$ i.e.
there exists a continuous $\Gamma$-equivariant map from 
$(X_{1},\rho)$ onto $(X_{2},\sigma)$. Generalizing the
corresponding results in [Wa1] and [Wa2] we obtain
the following.
\newline
\newline
{\bf Theorem 3 : }
{\it Let $X_{1},X_{2}$ be
nilmanifolds and 
$\rho ,\sigma$ be automorphism actions of a discrete abelian group
$\Gamma$
on $X_{1}$ and $X_{2}$
respectively. Suppose that 
$(X_{2},\sigma)$ is a factor of $(X_{1},\rho)$ and either 
$X_{1} = X_{2}$
or $X_{2}$ is a torus. 
Then there is a non-affine continuous $\Gamma$-equivariant map from 
$(X_{1},\rho )$ to $(X_{2}, \sigma )$  if and only if
$(X_{2},\sigma)$ is not ergodic.
 \/}
\section{Almost equivariant maps}
In this section  we give a necessary
and
sufficient condition for
 existence of a non-affine
almost equivariant map from
 $(X_{1},\rho )$ to $(X_{2},\sigma )$. 
Throughout this section for $i = 1,2$; $X_{i} = G_{i}/D_{i}$
will denote a nilmanifold, $\pi_{i}$ will denote the 
projection map from $G_{i}$ to $X_{i}$, $e_{i}$ will denote the
identity element of $G_{i}$ and $\bar{e_{i}}$ will denote the
image of $e_{i}$ in $X_{i}$ under the map $\pi_{i}$. It is known
 that in this case any homomorphism from $D_{1}$ to $D_{2}$ can
be extended to a continuous homomorphism from $G_{1}$ to $G_{2}$
 (cf. [Ma]). We will also use the following fact (cf. [AGH], pp. 54).
\newline
\newline
{\bf Proposition 2.1}  ( [AGH] ) {\bf : } 
{\it Let $X = G/D$ be a nilmanifold. For $a = (a_{1},\ldots ,a_{n})$ in
 ${\mathbb{R}}^{n}$ let $I_{a}$ denote the set defined by
$$ I_{a} = \{\ x\ |\ a_{i}\le x_{i}\le a_{i}+1\ \forall i = 1,\ldots ,n\
\}.$$
Then there exists an invertible linear map $T$ from ${\mathbb{R}}^{n}$
to $L(G)$ such that for all $a$ in ${\mathbb{R}}^{n}$, the set
$ {\rm exp}\circ T(I_{a})$ is a fundamental domain for $X = G/D$.
 \/}
\newline
\newline
The following proposition was proved in [Wa2]. (see also [AP]).
\newline
\newline
{\bf Proposition 2.2 : }
{\it Let $X_{1} = G_{1}/D_{1}, X_{2} = G_{2}/D_{2} $ be
nilmanifolds and $ f : X_{1}\rightarrow X_{2}$
 be a  continuous map. Let  $ F : G_{1}\rightarrow G_{2}$ be a lift of
$f$.
 Then there exist a $g_{0}\in{G_{2}}$, a continuous homomorphism
$\theta(f) : G_{1} \rightarrow G_{2}$ and a continuous map
$ P(f) : G_{1} \rightarrow G_{2}$ such that 
\medskip

a) $ P(f)(e_{1}) =  e_{2},\ P(f)(g\cdot\gamma) = P(f)(g)\ \ \forall g\in{\Gamma}$.

b)   $F(g) =  P(f)(g)\cdot g_{0}\cdot \theta(f)(g)\ \ \forall
g\in{\Gamma}$.

\medskip
\noindent
Moreover for a given $f$, the maps $\theta(f)$ and $P(f)$ are unique. \/}
\newline
\newline
Suppose $X_{1}$, $X_{2}$ are nilmanifolds and $\rho$, $\sigma$ are
automorphism actions of a discrete group $\Gamma$ on $X_{1}$ and
$X_{2}$ respectively. Let $\overline{\rho}$ and
$\overline{\sigma}$ denote the induced automorphism actions of
$\Gamma$
on $G_{1}$ and $G_{2}$ respectively. Then from the uniqueness part
of Proposition 2.2 it is follows that for any $\Gamma$-equivariant
continuous map $f$ from $(X_{1},\rho)$ to $(X_{2},\sigma)$, $P(f)$ is 
a $\Gamma$-equivariant
continuous map  from $(G_{1},{\overline{\rho}})$ 
to $(G_{2},{\overline{\sigma}})$.
Now we obtain the following.

\bigskip
\noindent
{\bf Lemma 2.1 :}
{\it Let $X_{1} = G_{1}/D_{1}, X_{2} = G_{2}/D_{2} $ be
nilmanifolds and $\rho $, $ \sigma $ be automorphism actions
of a discrete group $\Gamma $ on $X_{1}$ and $X_{2}$ respectively.
Then there exists a nonaffine continuous $\Gamma$-equivariant map from
$(X_{1},\rho )$ to
$(X_{2},\sigma )$ if and only if
there exists a nonzero continuous $\Gamma$-equivariant map $S$ from
$(X_{1},\rho )$ to $(L(G_{2}),\sigma_{e})$ such that
$ S(\bar{e_{1}}) = 0 $.
\/}

\bigskip
\noindent
{\bf Proof :}
Suppose there exists a nonaffine $\Gamma$-equivariant continuous map
$f$ from $(X_{1},\rho )$ to
$(X_{2},\sigma )$.
Let $P = P(f) : G_{1}\rightarrow G_{2}$ be as defined above. Since
$ P(e_{1}) =  e_{2}$ and $ P(g\cdot\gamma) = P(g)$ for all
$g\in{\Gamma}$,
 there exists a unique continuous map $Q$ from $X_{1}$ to
$G_{2}$ such that $Q(\overline{e_{1}}) = e_{2}$ and 
$P = Q \circ \pi_{1}$.
 It is easy to see that $Q$ is a
$\Gamma$-equivariant map from $(X_{1},\rho)$ to
$(G_{2},\overline{\sigma})$. Note that since $G_{2}$ is a connected 
simply connected nilpotent Lie group, the map
${\rm exp} : L(G_{2})\rightarrow G_{2}$ is a diffeomorphism.
Hence there is a unique map $S : X_{1}\rightarrow L(G_{2})$ such  
that $Q =  {\rm exp}\circ S$.
 Now $S(\overline{e_{1}}) = 0$ and
since exp is a $\Gamma$-equivariant map from $(L(G_{2}),\sigma_{e})$ to
$(G_{2},\overline{\sigma})$
 it is easy to see that $S$ is a $\Gamma$-equivariant map from
 $(X_{1},\rho)$ to  $(L(G_{2}),\sigma_{e})$. Since $f$ is a nonaffine
map,
$P(f) = {\rm exp}\circ S\circ\pi_{1}$ is non-constant i.e. $S$
is a nonzero map.

Now suppose there exists a non-zero  $\Gamma$-equivariant continuous
map $S$
from
$(X_{1},\rho )$ to $(L(G_{2}),\sigma_{e})$ such that
$S({\bar{e_{1}}}) = 0$. Define a map $f :X_{1}\rightarrow X_{2}$ by
$$ f(x) = \pi_{2}\circ {\rm exp}\circ S(x)\ \forall x\in{X_{1}}.$$   
It is easy to check that $f$ is a $\Gamma$-equivariant map from
$(X_{1},\rho )$ to $(X_{2},\sigma )$ and $P(f) =
{\rm exp}\circ S\circ\pi_{1}$. Since the map ${\rm exp}\circ S$ is
non-constant, so is $P$. Now from
the uniqueness part of
Proposition 2.2 it follows that $f$ is a nonaffine map.
\newline
\newline
{\bf Lemma 2.2 : }
{\it Let $X = G/D$ be a nilmanifold and $V$ be a finite dimensional vector
space over $\mathbb{R}$. Let $\Gamma$ be a discrete group and
$\rho $,
$\sigma $ be automorphism actions of $\Gamma$ on
$X$ and $V$ respectively. Then for any $\Gamma$-equivariant map
$S : X\rightarrow V$ there exists a finite index subgroup
$\Gamma_{0}\subset\Gamma$ such that
$S\circ\rho(\gamma) = S \ \ \forall \gamma\in{\Gamma_{0}}.$
 \/}
\newline
\newline
{\bf Proof :}
We define $A, A_{1},A_{2},\ldots \subset L(G)$ by
$$ A_{i} = \{ v\  |\  {\rm exp}(iv)\in{D}\}\ ,\ \ A = \cup A_{i}.$$
If $\pi: G\rightarrow G/D$ denotes the projection map then we 
define $B, B_{1},B_{2},\ldots \subset X$ by
$$ B_{i} = \pi\circ {\rm exp}(A_{i})\ ,\ \ B = \cup B_{i}.$$
From Proposition 2.1 it follows that each $B_{i}$ is a finite subset of 
$X$ and $B$ is dense in $X$. Also
 it is easy to see that 
 each $B_{i}$ is invariant under the action $\rho $. Therefore for any
 element $x$ in $B$, the $\rho$-orbit of $x$ is finite. 
Let $W$ denote the subspace of $V$ which consists of all elements
of $V$ whose $\sigma$-orbit is finite. Since $S$ is $\Gamma$-equivariant
and $B$ is a dense subset of $X$, it
follows that the image of $S$ is contained in $W$.
Now choose a basis $\{ w_{1},w_{2},\ldots w_{l}\}$ of $W$. Define
$\Gamma_{1},\Gamma_{2},\ldots ,\Gamma_{l}$ and $\Gamma_{0} $
by  
$$ \Gamma_{i} = \{ \gamma \in{\Gamma}\ |\ \sigma(\gamma )(w_{i})=
w_{i}\ \}
,\ \ \Gamma_{0} =  \cap \Gamma_{i} .$$
Since each $\Gamma_{i}\subset \Gamma $ is a subgroup of
finite index, so is $\Gamma_{0}$.
Since $\Gamma_{0}$ acts trivially on $W$ and image of $S$ is contained in
$W$, we conclude that $S$ is a $\Gamma_{0}$-invariant map.
\newline
\newline
{\bf Proof of Theorem 1 : }
Suppose there exists a finite index subgroup
$\Gamma_{0}\subset\Gamma$ and a nonaffine continuous
 map $f$ from $(X_{1},\rho )$ to
$(X_{2},\sigma )$  which is $\Gamma_{0}$-equivariant.
Then by Lemma 2.1 there exists a
nonzero continuous $\Gamma_{0}$-equivariant map $S$ from
$(X_{1},\rho )$ to $(L(G_{2}),\sigma_{e})$ such that
$ S(\bar{e_{1}}) = 0 $. 
Let $W$ denote the subspace of $L(G_{2})$ which consists of all elements
of $L(G_{2})$ whose $\sigma_{e}$-orbit is finite.
Then  from Lemma 2.2 it follows that the image of $S$ is
contained in $W$. Hence there exists a nonzero vector $v$ in
$L(G_{2})$ such that the $\sigma_{e}$-orbit of $v$ is finite.
To prove a) we choose a norm $||.||$ on $W$ and define
a function $p : W\mapsto {\mathbb{R}}$ by
$$ p(w) = {\rm Inf}\ \{ ||\sigma_{e}(\gamma )(w)||\ |\ \gamma\in{\Gamma}\
\}.$$
Since $\sigma_{e}$-orbit of any element of $W$ is finite, the map $q =
p\circ S$ is a non-constant
continuous $\Gamma$-invariant function from $X_{1}$ to ${\mathbb{R}}$. 

Now suppose both the conditions a) and b) are satisfied. Then $W$ is
a nonzero subspace of $L(G_{2})$ and there exists a finite index subgroup
$\Gamma_{0}\subset \Gamma$ such that the $\sigma_{e}$-action
of $\Gamma_{0}$ on $W$ is trivial. Let $q : X_{1}\mapsto {\mathbb{R}}$
be a non-constant continuous $\Gamma$-invariant function from $X_{1}$ to
${\mathbb{R}}$ and let $h : {\mathbb{R}}\mapsto W$ be a continuous
map such that the map $h\circ q$ is nonzero and
$h\circ q(\bar{e_{1}}) = 0$. Then $S = h\circ q$ is a nonzero continuous
$\Gamma_{0}$-equivariant map from
$(X_{1},\rho )$ to $(L(G_{2}),\sigma_{e})$ and
$ S(\bar{e_{1}}) = 0 $.
Applying Lemma 2.1 we see that  
there exists a nonaffine continuous $\Gamma_{0}$-equivariant map from
$(X_{1},\rho )$ to
$(X_{2},\sigma )$. 
\newline
\newline
Let $(X,d)$ be a metric space and $\rho$ be a continuous action of a group
$\Gamma$ on $X$. Then $(X,\rho)$ is said to be ${\it expansive \/}$
 if there exists $\epsilon > 0$ such that for any two distinct points
$x,y$ in $X$,
$$ {\rm Sup}\ \{\ d(\rho(\gamma)(x),\rho(\gamma)(y))\ |\ \gamma\in{\Gamma}\ \}
 \ge \epsilon .$$
Any such $\epsilon$ is called an expansive constant of $(X,\rho)$.
It is easy to check that the notion of expansiveness is independent
of the metric $d$.
\newline
\newline
Now as a corollary of Theorem 1 we obtain the following.
\newline
\newline
{\bf Corollary 2.1 : }
{\it Let $X_{1} = G_{1}/D_{1}, X_{2} = G_{2}/D_{2} $ be
nilmanifolds and $\rho $, $ \sigma $ be automorphism actions
of a discrete group $\Gamma $ on $X_{1}$ and $X_{2}$ respectively.
Suppose that either $(X_{1},\rho )$ is ergodic or $(X_{2},\sigma )$ 
is expansive.
Then every  continuous $\Gamma$-equivariant map from
$(X_{1},\rho )$ to
$(X_{2},\sigma )$ is an affine map. \/}
\newline
\newline
{\bf Proof : }
If  $(X_{1},\rho )$ is ergodic then there is no non-constant 
$\Gamma$-invariant continuous function from $(X_{1},\rho )$ to
${\mathbb{R}}$. Applying Theorem 1 we see that 
there exists no nonaffine continuous $\Gamma$-equivariant map from
$(X_{1},\rho )$ to
$(X_{2},\sigma )$ .

 Suppose that $(X_{2},\sigma )$
is expansive.
 Choose a metric $d$ on $X_{2}$ and an expansive constant $\epsilon > 0$
with respect to $d$. Define open sets $U\subset X_{2} $ and
$V\subset L(G_{2})$ by
$$ U = \{x\ |\ d(\bar{e_{2}},x) < \epsilon\ \}\ ,\ \ 
V = (\pi_{2}\circ {\rm exp})^{-1}(U).$$
We claim that for every nonzero vector $v$ in $L(G_{2})$, the
$\sigma_{e}$-orbit of $v$ is infinite. To see this choose any vector 
$v_{0}$   in $L(G_{2})$ such that the
$\sigma_{e}$-orbit of $v_{0}$ is finite.
Choose $\alpha > 0$ sufficiently small so that the
$\sigma_{e}$-orbit of $\alpha v_{0}$ is contained in $V$ and
 does not intersect the set ${\rm exp}^{-1}(D_{2}) - \{ 0\}$. Then 
the
$\sigma$-orbit of the element $x_{0} = \pi_{2}\circ {\rm exp}(\alpha v_{0})$
is contained in $U$. Since $\bar{e_{2}}$ is fixed by the action 
$\sigma$, it follows that $x_{0} = \bar{e_{2}}$ i.e.
$v_{0} = 0$. Now applying Theorem 1 we see that
every continuous $\Gamma$-equivariant map from
$(X_{1},\rho )$ to
$(X_{2},\sigma )$ is an affine map.
\section{Rigidity of toral automorphisms}
In this section we will consider automorphism actions of discrete groups
on tori.
 Suppose
$\rho $ is an automorphism action of a discrete group $\Gamma $ on
$T^{m}$. Then  $\widehat{\rho}$ 
will denote the
 automorphism action  of $\Gamma $ on $T^{m}$ defined by
$$ \widehat{\rho}(\gamma )(\chi ) = \chi \circ \rho(\gamma )
\ \ \forall \chi \in{{\widehat{T}}^{m}},\gamma \in{\Gamma}.$$
It is well known that $(T^{m},\rho )$ is ergodic if and only if
$\widehat{\rho}$ has no nontrivial finite orbit. Recall that if
$(T^{m},\rho )$ is
not ergodic then $F_{\rho}\subset \widehat{T}^{m}$  will denote
the subgroup consisting of all elements with finite $\widehat
{\rho}$-orbit and 
$\Gamma_{\rho} \subset \Gamma $ will denote the subgroup defined by
$$\Gamma_{\rho} = \{ \gamma\ |\ \chi\circ \rho(\gamma) = \chi \ \forall
\chi
\in{F_{\rho}}\ \}.$$
Since  $F_{\rho}$ is a finitely
generated group, it follows that  $\Gamma_{\rho} \subset \Gamma $ is
a 
subgroup of finite index.
\newline
\newline
{\bf Lemma 3.1 : }{\it Suppose $\Gamma$, $\rho$ and $\Gamma_{\rho}$
are as  above. Then there
exists  $\chi_{0}\in{\widehat{T}^{m}}$ and 
 $x_{0}\in{T^{m}}$ such that
$$\Gamma_{\rho}\  =\  \{ \gamma\ |\ \chi_{0}\circ \rho(\gamma)(x_{0}) =
 \chi_{0}(x_{0}) \}$$
 \/}
\newline
{\bf Proof :}
For any $\chi $ in $F_{\rho}$, let
 $\Gamma_{\chi }\subset\Gamma $ denote the stabilizer of $\chi$
under the $\Gamma$-action $\widehat{\rho}$.
We claim that for any $\chi_{1},\chi_{2}$ in $F_{\rho}$, there exists a
$\chi^{'}$
in $F_{\rho}$ such that $\Gamma_{\chi^{'} } =$
$ \Gamma_{\chi_{1} }\cap\Gamma_{\chi_{2} }$.
To see this, for $i = 1,2$ define $A_{i}\subset {\widehat{T}}^{m}$ by
$$ A_{i}\  =\  \{ \chi_{i}\circ\rho(\gamma )- \chi_{i}\ |\
\gamma\in{\Gamma}\ \}.$$
Since $\chi_{1},\chi_{2}$ are elements of ${F_{\rho}}$, both $A_{1}$ and
$A_{2}$ are
finite.
Choose
$n$ large enough so that $nA_{1}\cap A_{2} = \{ 0\}$. Define
$\chi^{'} = n\chi_{1} - \chi_{2}$. Clearly
 $ \Gamma_{\chi_{1} }\cap\Gamma_{\chi_{2} }$  
is contained in $\Gamma_{\chi^{'} }$.
On the other hand if $\gamma \in{\Gamma_{\chi^{'}}}$ then
$$n( \chi_{1}\circ\rho(\gamma )- \chi_{1}) =
  \chi_{2}\circ\rho(\gamma )-\chi_{2}.$$
Since $nA_{1}\cap A_{2} = \{ 0\}$, this implies that
$ \gamma \in{\Gamma_{\chi_{1}}\cap\Gamma_{\chi_{2}}}$.

Suppose  $ \chi_{1},\ldots ,\chi_{d}$
is a finite set of generators of $F_{\rho}$. From the above claim it
follows that
that  there exists a $\chi_{0}$ in $F_{\rho}$ such that
$$\Gamma_{\chi_{0}} = \Gamma_{\chi_{1}}\cap\cdots\cap\Gamma_{\chi_{d}}
 = \Gamma_{\rho}.$$
 Let $x_{0}\in{T^{m}}$ be any
element such that the cyclic subgroup generated by $x_{0}$
is dense in $T^{m}$. Then it is easy to see that
$$ \Gamma_{\rho}\  =\  \Gamma_{\chi_{0}}\  =\  \{ \gamma\ |\ \chi_{0}\circ
\rho(\gamma)(x_{0}) =
 \chi_{0}(x_{0}) \}.$$
\newline
{\bf Lemma 3.2 : }
{\it Let $\Gamma_{1}\subset \Gamma $ be a  subgroup of finite index
and $f$ be a $\Gamma_{1}$-invariant continuous map from $(T^{m},\rho)$ to
a
metric space $(Y,d)$. Then $f$ is $\Gamma_{\rho}$-invariant. \/}
\newline
\newline
{\bf Proof : }
First let us assume that $(Y,d) = {\mathbb{C}}$ with the usual metric.
Let $\widehat{f}:{\widehat{T}}^{m}\rightarrow {\mathbb{C}}$ be the Fourier
transform
of $f$. It is easy to check that $f$ is invariant under a
subgroup $\Gamma_{2}\subset \Gamma$ if and only if $\widehat{f}$ is
constant on each $\Gamma_{2}$-orbit under the action $\widehat{\rho}$.
Since 
$\Gamma_{\rho}$ acts trivially on $F_{\rho}$ under the action 
$\widehat{\rho}$, to
prove  
$\Gamma_{\rho}$-invariance of $f$ it is sufficient to show that 
$\widehat{f} = 0$ on $\widehat{T}^{m}- F_{\rho}$.

Since $\Gamma_{1}$ is a subgroup of finite index, for any $\phi$
in $\widehat{T}^{m}- F_{\rho}$, the $\Gamma_{1}$-orbit of $\phi$ is
infinite.
Since $f$ is  $\Gamma_{1}$-invariant, $\widehat{f}$ is constant on
the  $\Gamma_{1}$-orbit of $\phi$. Since
$\sum_{\phi}|\widehat{f}(\phi)|^{2} < \infty $, we conclude that
$\widehat{f}(\phi) = 0$.

Now let $(Y,d)$ be any arbitrary metric space and $f$ be a continuous
$\Gamma_{1}$-invariant function from $T^{m}$ to $Y$. If 
$\  C(Y,{\mathbb{C}})$ denotes the set of all continuous functions
from $Y$ to ${\mathbb{C}}$, then for each
$g$ in $ C(Y,{\mathbb{C}})$ the map $\  g\circ f$ is
$\Gamma_{1}$-invariant.
Since
 $\  C(Y,{\mathbb{C}})$ separates points of $Y$, from the previous
argument
it
follows
that $f$ is $\Gamma_{\rho}$-invariant.
\newline
\newline
{\bf Proof of Theorem 2 : } Suppose there exists a
non-zero
continuous $\Gamma$-equivariant map from $(T^{m},\rho )$ to
 $(X,\sigma)$. Then the condition a) follows from Corollary 2.1.
 Also from  Lemma 2.1 it follows
that there exists a 
non-zero     
continuous $\Gamma$-equivariant map $S$ from $(T^{m},\rho )$ to
 $(L(G),\sigma_{e})$. Applying  Lemma 2.2 and Lemma 3.2 we see 
that $S$ is
$\Gamma_{\rho}$-invariant. 
 This implies that the $\sigma_{e}$-action of $\Gamma_{\rho}$ on the image
of $S$ is
trivial.
Now the condition b) follows from the fact that $S$ is nonzero. 

Now suppose the conditions a) and b) are satisfied. Fix a finite subset 
$A = \{ \gamma_{1},\ldots ,\gamma_{d}\}$ of $\Gamma$ which contains
exactly one element of each right coset of $\Gamma_{\rho}$.
Let $W$ denote the
subspace of
$L(G)$ which is fixed by
 $\sigma_{e} (\gamma )$ for all $\gamma $ in $\Gamma_{\rho}$.
For any $\Gamma_{\rho}$ invariant map $h:T^{m}\rightarrow W$,
we define a map $h_{A}:T^{m}\rightarrow L(G)$ by
$$ h_{A} = \sum_{\gamma\in{A}} \sigma_{e}(\gamma^{-1})\circ h\circ
\rho(\gamma ).$$
Let $\gamma_{1}$ and $\gamma_{2} = \gamma_{0}\gamma_{1}$ be two
elements of $\Gamma $ belonging to the same right coset of
$\Gamma_{\rho}$.
Since $h$ is $\Gamma_{\rho}$-invariant and $\Gamma_{\rho}$-action on $W$
is trivial, it is easy to see that  
$$ \sigma_{e}(\gamma_{2}^{-1})\circ h\circ\rho(\gamma_{2})
= \sigma_{e}(\gamma_{1}^{-1})\circ h\circ\rho(\gamma_{1}).$$
Therefore if $B$ is another set containing exactly one element of
each coset of $\Gamma_{\rho}$ then $h_{A} = h_{B}$. Now it is easy to
verify
that for all $\gamma $ in $\Gamma $, 
$$h_{A}\circ\rho(\gamma ) = \sigma_{e}(\gamma)\circ h_{\gamma A}.$$
 Hence  $h_{A}$
is a $\Gamma$-equivariant map from $(T^{m},\rho )$ to
$(L(T^{m}),\sigma_{e} )$. We will show that for a suitable choice of $h$,
$h_{A}$ is nonzero and $h_{A}(e) = 0 $.

Let $\chi_{0}\in{{\widehat{T}}^{m}}$ and $x_{0}\in{T^{m}}$ be as in
Lemma 3.1.
Define $c_{0}, c_{1},\ldots ,c_{d}\in{S^{1}}$ by
$$c_{0} = 1,\ \ c_{i} = \chi_{0}\circ
\rho(\gamma_{i})(x_{0})\ \forall i = 1,\ldots ,d .$$
 Then $1,c_{1},\ldots,c_{d}$ are distinct. We choose 
a continuous map $g$ from $S^{1}$ to $W$ such that 
$$ g(c_{d}) \ne 0,\ \ g(c_{i}) = 0,\ \ i = 0,\ldots ,d-1 .$$
  Since the map $ g\circ \chi_{0} : T^{m}\rightarrow W$  is
$\Gamma_{\rho}$-invariant, from the previous argument it follows that the
map $S = (g\circ \chi_{0})_{A}$
 is a $\Gamma$-equivariant
map from $(T^{m},\rho )$ to $(L(T)^{n},D\sigma )$. 
Also it is easy to see that $S$ is nonzero and $S(e) = 0$.
 Now Theorem 2
follows from Lemma 2.1.
\newline
\newline
The following corollary generalizes earlier results of [AP]  and
[Wa1].
\newline
\newline  
{\bf Corollary 3.1 : }
{\it Let $A$ and $B$ be elements of $GL(m,{\mathbb{Z}})$ and
$GL(n,{\mathbb{Z}})$ respectively. Let $k_{A}$ be the smallest
positive integer $i$ such that $A^{i}$ has no eigenvalue which
is a proper root of unity. Then the following two  are
equivalent.

a) There exists a continuous nonaffine map 
$f : T^{m}\rightarrow T^{n}$ satisfying 

$\ \ f\circ A = B\circ f$.

b) 1 is an eigenvalue of $B^{k_{A}}$.
\/}
\newline
\newline
{\bf Proof : }
Let $\Gamma$ be the cyclic group, $\rho$ be the $\Gamma$-action on $T^{m}$
generated by $A$ and 
$\sigma$ be the $\Gamma$-action on $T^{n}$ generated by $B$.
 Then after suitable identifications we have,
$$ {\hat{T}}^{m} = {\mathbb{Z}}^{m},\ \ F_{\rho} = 
\{ z\in{{\mathbb{Z}}^{m}}\ |\ A^{i}(z) = z \ for\ some\ i\ \}.$$
It is easy to see that $A^{k_{A}}$ leaves $F_{\rho}$ invariant.
Since  no eigenvalue of $A^{k_{A}}$ is a proper root of unity, it
follows that $A^{k_{A}}$ leaves $F_{\rho}$ pointwise fixed. Suppose
$j$ is another positive integer such that $A^{j}$ 
leaves $F_{\rho}$ pointwise fixed. Then it easy to check that $j$ is
a multiple of ${k_{A}}$. Therefore $\Gamma_{\rho} = k_{A}{\mathbb{Z}}$.
 It is easy to see that the action $\sigma|_{\Gamma_{\rho}}$ has a
nonzero fixed point
in $L({\mathbb{R}}^{n})$ if and only if 1 is an eigenvalue of 
$B^{k_{A}}$. Now the given assertion follows from Theorem 2.
\section{Rigidity of factor maps}
In this section we will consider the case when $\Gamma$ is abelian
 and $X_{2}$ is a topological factor
of $X_{1}$ i.e. there exists a continuous $\Gamma$-equivariant map from
$X_{1}$ onto $X_{2}$.
We will need the following two results.
\newline
\newline
{\bf Theorem 4.1} ( see [Be], Theorem 5.1){\bf :} 
{\it Let $\Gamma$ be a discrete
abelian group, $T^{n}$ be
the $n$-torus  and $\rho$ be an
ergodic
automorphism action of $\Gamma$ on $T^{n}$. Then there exists an element
$\gamma_{0}$ of $\Gamma$ such that $\rho(\gamma_{0})$ is an 
ergodic automorphism. \/}
\newline
\newline
{\bf Theorem 4.2 }( see [Pa]){\bf :}
{\it Let $X =G/D$ be
a nilmanifold and $\theta$ be an 
automorphism  of $X$ such that $\theta$ induces an ergodic
automorphism on the torus $G/ [G,G]\circ D$. Then
 $\theta$ is an ergodic automorphism of $X$. \/}
\newline
\newline
If $X = G/D$ is a nilmanifold then by $X^{0}$ we denote the torus
$G/[G,G]\cdot D$ and by $\pi^{0}$ we denote the projection map
from $G$ onto $X^{0}$. If $\rho$ is an automorphism action of a discrete
group $\Gamma$ on $X$ then $\rho^{0}$ will denote the 
automorphism action of $\Gamma$ on $X^{0}$ induced by $\rho$.
\newline
\newline
We note the following simple consequence of Theorem 4.1 and Theorem 4.2.
\newline
\newline
{\bf Proposition 4.1 :}
{\it Let $\Gamma$ be a discrete abelian group, X = G/D be a nilmanifold
 and $\rho$ be an automorphism action of $\Gamma$ on $X$. Then
$(X,\rho)$ is ergodic if and only if
$(X^{0},\rho^{0})$ is ergodic. \/}
\newline
\newline
{\bf Proof : }
Let $q : X\rightarrow X^{0}$ denote the projection map. Then it is
easy to check that $q$ is a measure preserving $\Gamma$-equivariant
map from $(X,\rho)$ to $(X^{0},\rho^{0})$. Therefore ergodicity of
$(X,\rho)$ implies  ergodicity of $(X^{0},\rho^{0})$. On the other
hand if $(X^{0},\rho^{0})$ is ergodic then by Theorem 4.1 there exists
a $\gamma$ in $\Gamma$ such that $\rho^{0}(\gamma )$ is an ergodic
automorphism of $X^{0}$. Applying Theorem 4.2 we see that 
$\rho (\gamma)$  is an ergodic
automorphism of $X$ i.e. $(X,\rho)$ is ergodic.
\newline
\newline
If $V$ is a finite dimensional vector space over a field $K$ then
by $V^{*}$ we denote the dual of $V$. If $\Gamma$ is a discrete group 
and $\rho : \Gamma \rightarrow GL(V)$ is an automorphism action of
$\Gamma$ on $V$ then by $\rho^{*}$ we denote the action of
$\Gamma$ on $V^{*}$ defined by
$$ \rho^{*}(\gamma)(q)(v) = q(\rho^{*}(\gamma)^{-1}v)
\ \ \forall q\in{V^{*}},v\in{V}.$$
\newline
{\bf Proposition 4.2 :}
{\it Let $\Gamma$ be an abelian group and $V$ be a finite dimensional
vector
space over $\mathbb{R}$. Let $\rho :\Gamma\rightarrow GL(V)$ be an 
automorphism action of $\Gamma$ on $V$ such that the
induced $\Gamma$-action on the dual of $V$ has a nontrivial fixed
point. Then $\rho$ has nontrivial fixed point in $V$. \/} 
\newline
\newline
{\bf Proof :}
By passing to the complexification we see that it is enough to prove
the analogous statement when $V$ is a finite dimensional vector space
over ${\mathbb{C}}$. In that case after suitable
identifications
we can assume that $V = V^{*} = {\mathbb{C}}^{n}$, 
$\rho :\Gamma\rightarrow GL(n,{\mathbb{C}})$ is a homomorphism and
$\rho^{*} :\Gamma\rightarrow GL(n,{\mathbb{C}})$ is the homomorphism defined
by $\rho^{*}(\gamma) = \rho(\gamma^{-1})^{T}$. Let us consider the special
case when with respect to some basis in ${\mathbb{C}}^{n}$ each
$\rho(\gamma)$ is given by an upper triangular matrix
with equal diagonal entries. In this case it is easy to verify that
$\rho$ or $\rho{*}$ has a nonzero fixed vector in ${\mathbb{C}}^{n}$
if and only if for any $\gamma$ in $\Gamma$ all the diagonal entries
of $\rho(\gamma)$ are equal to 1. To prove the general case we note that
since $\Gamma$ is abelian, there exist subspaces 
$V_{1},V_{2},\ldots ,V_{k}$ of ${\mathbb{C}}^{n}$ and homomorphisms
$\rho_{i} :\Gamma \rightarrow GL(V_{i})$; $i= 1,\ldots ,k$ such that
${\mathbb{C}}^{n} = V_{1}\oplus \cdots \oplus V_{k}$,
$\rho = \rho_{1}\oplus \cdots \oplus \rho_{k}$ and
each $\rho_{i}$ satisfies the above condition (cf. [Ja], pp. 
134).
\newline
\newline
{\bf Proposition 4.3 :}
{\it Let $\sigma$ be an automorphism action of a discrete abelian group
$\Gamma$ on a torus $T^{n}$. Then $(T^{n},\sigma )$ is ergodic if and only 
if there is no nonzero element in $L(T^{n})$ with finite 
$\sigma_{e}$-orbit. \/}
\newline
\newline
{\bf Proof : }
Since
$T^{n} = {\mathbb{R}}^{n}/{\mathbb{Z}}^{n}$, 
$L(T^{n})$ can be identified with ${{\mathbb{R}}^{n}}$.
 Also $\sigma_{e}$ can be   
realised as a homomorphism from $\Gamma$ to $GL(n,{\mathbb{Z}})$,
the dual action
$\sigma_{e}^{*}$ can be
realised as the homomorphism from $\Gamma$ to $GL(n,{\mathbb{Z}})$
which takes $\gamma$ to $ \sigma_{e}(\gamma^{-1})^{T}$
and ${\widehat{\sigma}}$
can be identified with $\sigma_{e}^{*}|_{{\mathbb{Z}}^{n}}$.
 Suppose $(T^{n},\sigma )$ is ergodic. 
Let $\Gamma_{0}\subset\Gamma$ be a subgroup of finite
index. Then no nonzero
element
of ${\mathbb{Z}}^{n}$ is fixed by $\Gamma_{0}$ under the action
$\sigma_{e}^{*}$. Since
$\sigma_{e}^{*}(\gamma)\in{GL(n,{\mathbb{Z}})}$ for all $\gamma$,  
this implies that no nonzero element
of ${\mathbb{R}}^{n}$ is fixed by $\Gamma_{0}$ under the action
$\sigma_{e}^{*}$. Applying Proposition 4.1 we see that
no nonzero element
of ${\mathbb{R}}^{n}$ is fixed by $\Gamma_{0}$ under the action
$\sigma_{e}$.
Now suppose $(T^{n},\sigma )$ is not ergodic. Then there exists a finite
index subgroup $\Gamma_{0}\subset\Gamma$ and a nonzero point $z$ in
${\mathbb{Z}}^{n}$ such
that  $z$ is fixed by $\Gamma_{0}$ under the action $\sigma_{e}^{*}$.
 Now from Proposition 4.1 we conclude that
there exists a nonzero element
in ${\mathbb{R}}^{n}$ which is fixed by $\Gamma_{0}$ under the action
$\sigma_{e}$.  
\newline
\newline
{\bf Proof of Theorem 3 : }
 Suppose $(X_{2},\sigma)$ is not ergodic. By our
assumption there exists a continuous $\Gamma$-equivariant map $f$ from
$(X_{1},\rho)$ onto $(X_{2},\sigma)$. If $f$ is nonaffine then there
is nothing to prove. Therefore we may assume that there exists a
$g_{0}\in{G_{2}}$ and a continuous homomorphism   
$\theta : G_{1}\rightarrow G_{2}$ such that
$f(gD_{1}) = g_{0}\theta(g)D_{2}$ for all $g$ in $G_{1}$.
Since $f$ is surjective and $\Gamma$-equivariant, so is $\theta$.
Let $\theta_{0}$ denote the homomorphism from $X_{1}^{0}$ to
$X_{2}^{0}$
 induced by $\theta$.
Then  $\theta_{0}$ is surjective.
 Since $(X_{2},\sigma)$ is not ergodic, from Proposition 4.1 it follows
that $(X_{2}^{0},\sigma^{0} )$ is not ergodic.
Let  $\phi$ be an element of $X_{2}^{0}$ such that
$\widehat{\sigma^{0}}$-orbit
of $\phi$ is finite. Since $\theta^{0}$ is $\Gamma$-equivariant, it
follows that $\widehat{\rho^{0}}$-orbit of $\phi\circ\theta^{0}$ is also
finite,
which implies that $(X_{1}^{0},\rho_{0} )$ is not ergodic.
Also for any $\gamma$ in $\Gamma_{\rho}$,
$$\phi\circ\sigma^{0}(\gamma)\circ\theta^{0} =
  \phi\circ\theta^{0}\circ\rho^{0}(\gamma) = \phi\circ\theta^{0}.$$
 Since $\theta^{0}$
is a surjective map, this implies that $\phi\circ\sigma^{0}(\gamma) = $
$\phi$ for all $\gamma$ in $\Gamma_{\rho}$. Let $\pi_{2}^{0}$ denote
the projection map from $G_{2}$ onto $X_{2}^{0}$ and let $q$ denote
the map $\phi\circ\pi_{2}^{0}\circ {\rm exp}$. Then
$dq : L(G_{2})\rightarrow {\mathbb{R}}$ is an element of the dual
of $ L(G_{2})$ such that $dq\circ\sigma_{e}(\gamma) = $   
$dq$ for all $\gamma$ in $\Gamma_{\rho}$.
 Now from Proposition 4.2 it follows that there exists a nonzero point
in
$ L(G_{2})$ which is fixed by $\Gamma_{\rho}$ under the action
$\sigma_{e}$. Applying Theorem 2 we see that
there exists a continuous nonaffine $\Gamma$-equivariant map $h$ from  
$(X^{0}_{1},\rho^{0})$ to
$(X_{2},\sigma )$. If $\pi_{1}^{0}$ denotes the projection map from
$X_{1}$
to $X^{0}_{1}$ then it is easy to see that $h\circ \pi_{1}^{0}$ is a
continuous nonaffine $\Gamma$-equivariant map $h$ from
$(X_{1},\rho)$ to
$(X_{2},\sigma )$.
 
Now suppose  $(X_{2},\sigma )$ is ergodic. Since by our assumption
either $(X_{1},\rho) = (X_{2},\sigma)$
or $X_{2}$ is a torus from Proposition 4.3 it follows that
 either  $(X_{1},\rho )$ is ergodic or there is no non-zero
 element in $L(G_{2})$ whose $\sigma_{e}$-orbit is finite.
 Applying Theorem 1 we conclude that
 every continuous $\Gamma$-equivariant map from $(X_{1},\rho)$ to
$(X_{2},\sigma )$ is an affine map.
\newline
\newline
The following examples show that Theorem 3 does not hold if any of
the assumptions as in the hypothesis is dropped.
\newline
\newline
{\bf Example 1 :}
Let $\Gamma$ be the cyclic group and $\rho ,\sigma$ be the
automorphism actions of $\Gamma$
on ${\mathbb{R}}/{\mathbb{Z}}$ generated by the identity automorphism and
the automorphism
$z\rightarrow -z $ respectively. Then it is easy to see that in this
case 
$\Gamma_{\rho} = \Gamma$
and no nonzero element of $L({\mathbb{R}})$ is fixed by $\Gamma_{\rho}$
under the action
$\sigma_{e}$. Now applying Theorem 2 we conclude that
there is no
 nonaffine continuous $\Gamma$-equivariant map from 
$(S^{1},\rho )$ to $(S^{1},\sigma )$. Note that in this case $\Gamma$ is
abelian and
neither of the two actions is ergodic.
\newline
\newline
{\bf Example 2 :}
Fix $n\ge 3$ and define a  subgroup $\Gamma$ of $GL(n,\mathbb{Z})$
 by
$$ \Gamma =  \left\{  
\left(
\begin{array}{cc}
 A & b  \\
 0 & 1  \\
\end{array}
\right) 
\ |\  A\in{GL(n-1,\mathbb{Z})}, b\in{{\mathbb{Z}}^{n-1}}\ 
 \right\}  $$
Let $\rho$ denote the natural action of $\Gamma$ on 
${\mathbb{R}}^{n}/{\mathbb{Z}}^{n}$. Then
it is easy to see that for any $x = (x_{1},\ldots ,x_{n})$ 
in $L({\mathbb{R}}^{n})$, the
$\rho_{e}$-orbit of $x$ is 
unbounded. Applying Theorem 1 we see that  
there is no 
 nonaffine continuous $\Gamma$-equivariant map from
$(T^{n},\rho)$ to $(T^{n},\rho)$. Note that in this case $(T^{n},\rho )$
is 
not ergodic since the vector $x_{0} = (0,\ldots ,0,1)$
is fixed by the dual action $\rho^{*}$.
\newline
\newline
{\bf Example 3 : }
Suppose $X = G/D$ where $G$ and $D$ are defined by
$$ G =  \left\{
{\scriptsize 
\left(
\begin{array}{ccc} 
 1 & x & z \\
 0 & 1 & y \\
 0 & 0 & 1 \\
\end{array}
\right)
}
\ |\  x,y,z \in{\mathbb{R}}\ \right\}\ ,\ 
D = \left\{ 
{\scriptsize  
\left(
\begin{array}{ccc}
 1 & p & r \\
 0 & 1 & q \\
 0 & 0 & 1 \\
\end{array}
\right)
}
\ |\  p,q,r \in{\mathbb{Z}}\ \right\}$$
 Let $A$ be an ergodic automorphism
of $G/D$. If $G_{0}$ denotes the center of $G$ then it is easy to see
that $G_{0}/G_{0}\cap D $ is isomorphic to $S^{1}$. Hence replacing
$A$ by $A^{2}$ if necessary we may assume that $A$ acts trivially
on $G_{0}$.  
 Define a nilmanifold
$X_{1}$ and an automorphism $A_{1}$ of $X_{1}$ by
$X_{1} = X \times S^{1} , A_{1} = A\times Id$. Let $\rho_{1}$
and $\rho$ denote the automorphism actions of $\mathbb{Z}$ on
$X_{1}$ and $X$ generated by $A_{1}$   
and $A$ respectively. Then $(X ,\rho )$ is a factor of
$(X_{1},\rho_{1})$. Let $\pi : X_{1}\rightarrow S^{1}$ be the
the projection map and
$h : S^{1}\rightarrow L(G_{0})$ be  any nonzero map such that
$h(e) = 0$. Then $h\circ \pi $ is a nonzero $\Gamma$-equivariant map
from $(X_{1},\rho_{1})$ to $(L(G),\rho_{e})$ such that 
$ h\circ \pi (e) = 0$. Applying Lemma 2.1
we see that
there exists a nonaffine continuous $\Gamma$-equivariant map
from $(X_{1},\rho_{1})$ to $(X,\rho)$.
\newpage
{ \Large \bf References : }
\newline
\newline
[ A-G-H ]\ \ L. Auslander, L. Green and F. Hahn (1963) {\it Flows on
homogeneous
spaces. \/} Princeton University Press.
\newline
[ A-P ]\ \ R.L Adler and R. Palais (1965) {\it Homeomorphic conjugacy
of automorphisms on the torus. \/} Proc. Amer. Math Soc.
{\bf 16}: 1222-1225.
\newline
[ Be ]\ \ D. Berend (1985) {\it Ergodic semigroups of epimorphisms. \/}
Trans. Amer. Math. Soc.
{\bf 289}: 393-407.
\newline    
[ Ja ]\ \ N. Jacobson (1953) {\it Lectures in abstract algebra,\/} Vol 2.
Van Nostrand.
\newline   
 [ Ma ]\ \ A. Mal'cev (1962) {\it On a class of homogeneous spaces. \/}
Amer. Math. Soc. Transl. Ser. (2) {\bf 39}: 276-307.
\newline
[ Pa ]\ \ W. Parry (1970) {\it Dynamical systems on nilmanifolds. \/}  
 Bull. London. Math. Soc.
{\bf 2}: 37-40.
\newline
[ Ra ]\ \ M.S Raghunathan (1972) {\it Discrete subgroups of Lie groups.
\/}
 Springer Verlag.
\newline
[ Wa1 ]\ \  P. Walters (1968) {\it Topological conjugacy of affine
transformations of tori. \/}
Trans. Amer. Math. Soc. {\bf 131}: 40-50.
\newline
[ Wa2 ]\ \  P. Walters (1970) {\it Conjugacy properties of affine
transformations on nilmanifolds. \/}
Math. Systems Theory  {\bf 4}: 326-337.
\newline
\newline
\newline
 Address : School of mathematics.\newline
Tata Institute of Fundamental Research.\newline
Homi Bhabha road.\newline
Mumbai - 400005, India.
\newline
\newline
 e-mail : siddhart@math.tifr.res.in

\end{document}